\documentclass[a4paper,10pt,leqno]{article}

\usepackage{amssymb,latexsym}
\usepackage{mltex}
\usepackage{amsmath}
 \usepackage{lscape}
\usepackage{enumerate}
\usepackage{amsthm}
\usepackage{hyperref}

\newtheorem{thm}{Theorem}

\newtheorem{lem}{Lemma}[section]
\newtheorem{prop}[lem]{Proposition}
\newtheorem{rem}{Remark}
\newtheorem{defn}[lem]{Definition}

\newtheorem*{thrm}{Theorem}

\newcommand{\Ekt}{\mathbb{E}(\kappa,\tau)}

\newcommand{\trace}{\mathrm{tr\,}}

\newcommand{\lto}{\ensuremath{\longrightarrow}}
\newcommand{\C}{\mathbb{C}}
\newcommand{\Dt}{\dfrac{\partial }{\partial t}}
\newcommand{\dt}{\frac{\partial }{\partial t}}
\newcommand{\Ss}{\mathbb{S}}
\newcommand{\HH}{\mathbb{H}}

\newcommand{\R}{\mathbb{R}}

\newcommand{\M}{\mathbb{M}}

\newcommand{\pre}{\Re e}

\newcommand{\base}{\{e_1,\ldots ,e_n\}}

\newcommand{\function}[5]
{\begin{eqnarray*}\begin{array}{r@{}ccl}
 #1\;\colon\;  & #2 &\lto & #3 \\[.05cm]
  & #4 &\longmapsto  & #5
\end{array}\end{eqnarray*}
}
\newcommand{\beqt}{\begin{equation}}  \newcommand{\eeqt}{\end{equation}}
\newcommand{\bal}{\begin{align}}      \newcommand{\eal}{\end{align}}
\newcommand{\ba}{\begin{array}}      \newcommand{\ea}{\end{array}}
\newcommand{\bc}{\begin{center}}     \newcommand{\ec}{\end{center}}
\newcommand{\be}{\begin{enumerate}}  \newcommand{\ee}{\end{enumerate}}
\newcommand{\beq}{\begin{eqnarray}}  \newcommand{\eeq}{\end{eqnarray}}
\newcommand{\beQ}{\begin{eqnarray*}} \newcommand{\eeQ}{\end{eqnarray*}}
\newcommand{\bi}{\begin{itemize}}    \newcommand{\ei}{\end{itemize}}
\newcommand{\bt}{\begin{tabular}}    \newcommand{\et}{\end{tabular}}
\newcommand{\finpreuve}{\hfill\square\\}

\def\pf{\noindent{\textit {Proof :} }}

\title{Isometric Immersions of Hypersurfaces in 4-dimensional Manifolds via Spinors}
\author{Marie-Am\'elie Lawn\;\, and Julien Roth}
\date{}
\begin{document}
\maketitle

\begin{abstract}
We give a spinorial characterization of isometrically immersed hypersurfaces into 4-dimensional space forms and product spaces $\M^3(\kappa)\times\R$, in terms of the existence of particular spinor fields, called generalized Killing spinors or equivalently solutions of a Dirac equation. This generalizes to higher dimensions several recent results for surfaces by T.\,Friedrich, B.\,Morel and the two authors. The main argument is the interpretation of the energy-momentum tensor of a generalized Killing spinor as the second fundamental form, possibly up to a tensor depending on the ambient space. As an application, we deduce some non-existence results for isometric immersions into the 4-dimensional Euclidean space
\end{abstract}
{\it keywords:} Dirac Operator, Generalized Killing Spinors, Isometric Immersions, Gauss and Codazzi-Mainardi Equations, Energy-Momentum Tensor.

\date{}
\maketitle\pagenumbering{arabic}
\section{Introduction}
A classical problem in Riemannian geometry is to know when a
Riemannian manifold $(M^n,g)$ can be
isometrically immersed into a fixed Riemannian manifold $(\bar{M}^{n+p},\bar{g})$. In this paper, we will focus on the case of hypersurfaces, that is $p=1$.\\
\indent The case of space forms $\R^{n+1}$, $\Ss^{n+1}$ and
$\HH^{n+1}$ is well-known. The Gauss and Codazzi-Mainardi equations
are necessary and sufficient conditions. Recently, B. Daniel
(\cite{Dan})
gave an analogous characterization for hypersurfaces in the product spaces $\Ss^n\times\R$ and $\HH^n\times\R$.\\
\indent In low dimensions, namely for surfaces, another necessary and
sufficient condition is now well-known, namely the existence of a
special spinor field called {\it generalized Killing spinor field}
(see \cite{Fr,Mo,Roth4,La,LR}). Note that this condition is not
restrictive since any oriented surface is also spin. This approach
was first used by T. Friedrich (\cite{Fr}) for surfaces in $\R^3$
and then extended to other $3$-dimensional Riemannian manifolds by
(\cite{Mo,Roth4}).\\
 \indent More generally, the restriction
$\varphi$ of a parallel spinor field on $\R^{n+1}$ to an oriented
Riemannian hypersurface $M^{n}$ is a solution of a generalized
Killing equation \beqt\label{killing} \nabla^{\Sigma
M}_X\varphi=-\frac{1}{2}\gamma^M(A(X))\varphi, \eeqt where
$\gamma^M$ and $\nabla^{\Sigma M}$ are respectively the Clifford
multiplication and the spin connection on $M^{n}$, and $A$ is the
Weingarten tensor of the immersion. Conversely, Friedrich proves in
\cite{Fr} that, in the two dimensional case, if there exists a
generalized Killing spinor field satisfying equation
\eqref{killing}, where $A$ is an arbitrary field of symmetric
endomorphisms of $TM$, then $A$ satisfies the Codazzi-Mainardi and
Gauss equations of hypersurface theory and is consequently the
Weingarten tensor of a local isometric immersion of $M$ into $\R^3$.
Moreover, in this case, the solution $\varphi$ of the generalized
Killing equation is equivalently a solution of the Dirac equation
\beqt\label{dirac} D\varphi=H\varphi, \eeqt
 where $|\varphi|$ is constant and $H$ is a real-valued function.\\
\indent One feature of those spinor representations is that
fundamental topological informations can be read off more easily
from the spinorial data (see for example \cite{KS}).\\
\indent The question of a spinorial characterization of
3-dimensional manifolds as hypersurfaces into a given 4-dimensional
manifold is also of special interest since, again, any oriented
3-dimensional manifold is spin. The case of hypersurfaces of the
4-dimensional Euclidean space has been treated by Morel in
\cite{Mo}, when $A$ is a Codazzi tensor. Here, we extend Morel's
result to other 4-dimensional space forms and product spaces, that
is $\Ss^4$,
$\HH^4$ (see Theorem \ref{thmS}),  $\Ss^3\times\R$ and $\HH^3\times\R$ (see Theorem \ref{thmSR1}).\\
\indent The techniques we use in this article are different
from those in Friedrich and Morel's approach. The main difference is
that unlike in the 2-dimensional case, the spinor bundle of a
3-dimensional manifold does not decompose into subbundles of
positive and negative half-spinors. In this case, the condition for
an isometric immersion is the existence of two particular spinor
fields on the manifold instead of one as in the case of surfaces.
Moreover, we prove the equivalence between the generalized Killing
equation and the Dirac equation for
spinor fields of constant norm in the above four cases.\\
\indent The last paragraph is devoted to an application. We prove in
a straightforward way using our results and the existence of special
spinors on certain three-dimensional $\eta$-Einstein manifolds that
they cannot be isometrically immersed into the Euclidean space
$\mathbb{R}^4$.
\section{Preliminaries}
\subsection{Hypersurfaces and induced spin structures}
\label{sect22}
We begin by preliminaries on hypersurfaces and induced spin structures. The reader can refer to \cite{LM,Fr3,BHMM} for basic facts about spin geometry and \cite{Ba2,Mo3,HMo} for the spin geometry of hypersurfaces.\\
\indent
Let $(N^{n+1},g)$ be a Riemannian spin manifold and $\Sigma N$ its spinor bundle.
We denote by $\nabla$ the Levi-Civita connection on $TN$, and $\nabla^{\Sigma N}$ the spin connection on $\Sigma N$.
The Clifford multiplication will be denoted by $\gamma$ and $\left\langle .,.\right\rangle $ is the natural Hermitian product on $\Sigma N$,
compatible with $\nabla$ and $\gamma$. Finally, we denote by $D$ the Dirac operator on $N$ locally given by $D=\sum_{i=1}^n\gamma(e_i)\nabla_{e_i}$, where $\{e_1,\cdots,e_{n+1}\}$ is an orthonormal frame of $TN$.\\
\indent
Now let $M$ be an orientable hypersurface of $N$. Since the normal bundle is trivial, the hypersurface $M$ is also spin. Indeed, the existence of a normal unit vector field $\nu$ globally defined on $M$ induces a spin structure from that on $N$.\\
\indent
Then we can consider the intrinsic spinor bundle of $M$ denoted by
$\Sigma M$. We denote  respectively by $\nabla^{\Sigma M}$,
$\gamma^{M}$ and $D^{M}$, the Levi-Civita connection, the Clifford multiplication and the intrinsic Dirac operator on $M$. We can also define an extrinsic spinor bundle on $M$ by ${\bf{S}}:=\Sigma N_{|M}$. Then we recall the identification between these two spinor bundles ({\it cf} \cite{HMo}, \cite{Mo3} or \cite{Ba2} for instance):
\beqt\label{identifications}
{\bf{S}}\equiv \left\lbrace
\begin{array}{ll}
 \Sigma M&\text{if $n$ is even}  \\
\Sigma M\oplus\Sigma M & \text{if $n$ is odd}.
\end{array} \right.
\eeqt \noindent The interest of this identification is that we can
use restrictions of ambient spinors to study the intrinsic Dirac
operator of $M$. Indeed, we can define an extrinsic connection
$\nabla^{{\bf{S}}}$ and a Clifford multiplication
$\gamma^{{\bf{S}}}$ on ${\bf{S}}$ by \beqt\label{gaussspin}
\nabla^{{\bf{S}}}=\nabla+\frac{1}{2}\gamma(\nu)\gamma(A), \eeqt
\beqt \gamma^{{\bf{S}}}=\gamma(\nu)\gamma, \eeqt \noindent where
$\nu$ is the exterior normal unit vector field and $A$ the
associated Weingarten operator. By the previous identification given
by (\ref{identifications}), we can also identify connections and
Clifford multiplications. \beqt \nabla^{{\bf{S}}}\equiv \left\lbrace
\begin{array}{ll}
 \nabla^{\Sigma M}& \text{if $n$ is even,} \\
\nabla^{\Sigma M}\oplus \nabla^{\Sigma N}& \text{if $n$ is odd,}
\end{array} \right.
\eeqt
\beqt
\gamma^{{\bf{S}}}\equiv \left\lbrace \begin{array}{ll}
 \gamma^{M}& \text{if $n$ is even,} \\
\gamma^{M}\oplus -\gamma^{M}& \text{if $n$ is odd.}
\end{array} \right.
\eeqt
\noindent Then, we can consider the following extrinsic Dirac operator on $M$, acting on sections of ${{\bf{S}}}$, denoted by ${\bf{D}}$ and given locally by
\beqt
{\bf{D}}=\sum_{i=1}^{n}\gamma^{{\bf{S}}}(e_i)\nabla^{{\bf{S}}}_{e_i},
\eeqt
\noindent where $\base$ is an orthonormal local frame of
$TM$. Then, by (\ref{gaussspin}), we have
\beqt
{\bf{D}}=\frac{n}{2}{H}-\gamma(\nu)\sum_{i=1}^{n}\gamma(e_i)\nabla_{e_i},
\eeqt
\noindent that is, for any $\psi\in\Gamma({\bf{S}})$
\beqt\label{Diracbord}
{\bf{D}}\psi:=\frac{n}{2}{H}\psi-\gamma(\nu)D\psi-\nabla_{\nu}\psi.
\eeqt
\begin{rem}
In the sequel, when we are only considering 3-dimensional manifolds, we will denote for the sake of simplicity the Clifford multiplication by a dot.
\end{rem}

We have all the spinorial ingredients, and now, we will give some reminders about surfaces into product spaces.
\subsection{Basic facts about product spaces}
In this section, we recall some basic facts on the product spaces $\M^n(\kappa)\times\R$ and their hypersurfaces. More details can be found in \cite{Dan} for instance. In the sequel, we will denote by $\M^n(\kappa)$ the $n$-dimensional simply connected space form of constant sectional curvature $\kappa$. That is,
$$\M^n(\kappa)=\left\{
\begin{array}{ll}
\Ss^n(\kappa)&\text{if}\ \kappa>0\\
\R^n&\text{if}\ \kappa=0\\
\HH^n(\kappa)&\text{if}\ \kappa<0.
\end{array}
\right.
$$
We denote by $\overline{\nabla}$ and $\overline{R}$ the Levi-Civita connection and the curvature tensor of $\M^n(\kappa)\times\R$. Finally, let $\dt$ be the unit vector field giving the orientation of $\R$ in the product $\M^n(\kappa)\times\R$.\\
Now, let $M$ be an orientable hypersurface of $\M^n(\kappa)\times\R$ and $\nu$ its unit normal vector. Let $T$ be the projection of the vector $\dt$ on the tangent bundle $TM$. Moreover, we consider the function $f$ defined by:
$$f:=\Big<\nu,\Dt\Big>.$$
It is clear that
$$\Dt=T+f\nu.$$
Since $\dt$ is a unit vector field, we have:
$$||T||^2+f^2=1.$$
Let's compute the curvature tensor of $\M^n(\kappa)\times\R$ for tangent vectors to $M$.
\begin{prop}{\cite{Dan,thRoth}}
For all $X,Y,Z,W\in\Gamma(TM)$, we have:
\begin{eqnarray*}
\big<\overline{R}(X,Y)Z,W\big>&=&\kappa\big(\langle X,Z\rangle \langle Y,W\rangle -\langle Y,Z\rangle \langle X,W\rangle \\
&&-\langle Y,T\rangle \langle W,T\rangle \langle X,Z\rangle -\langle X,T\rangle \langle Z,T\rangle \langle Y,W\rangle \\
&&+\langle X,T\rangle \langle W,T\rangle \langle Y,Z\rangle +\langle Y,T\rangle \langle Z,T\rangle \langle X,W\rangle \big),
\end{eqnarray*}
and
\begin{eqnarray*}
\big<\overline{R}(X,Y)\nu,Z\big>&=&\kappa f\big(\langle X,Z\rangle \langle Y,T\rangle -\langle Y,Z\rangle \langle X,T\rangle \big).
\end{eqnarray*}
\end{prop}
\noindent
The fact that $\Dt$ is parallel implies the following two identities
\begin{prop}{\cite{Dan,thRoth}}
For $X\in\Gamma(TM)$, we have
\begin{eqnarray}\label{FirstequProductsp}
\nabla_XT=fA(X),
\end{eqnarray}
and
\begin{eqnarray}\label{SecondequProductsp}
df(X)=-\langle A(X),T\rangle.
\end{eqnarray}
\end{prop}
\noindent
{\it Proof:} We know that $\overline{\nabla}_X\dt=0$ and $\dt=T+f\nu$,
so
\begin{eqnarray*}
0&=&\overline{\nabla}_XT+df(X)\nu+f\overline{\nabla}_X\nu\\
&=&\nabla_XT+\langle A(X),T\rangle \nu+df(X)\nu-fA(X).
\end{eqnarray*}
Now, it is sufficient to consider the normal and tangential parts to obtain the above identities.
$\finpreuve$
\begin{defn}[Compatibility Equations]\label{comp}
We sA(Y) that $(M,\langle.,.\rangle,A,T,f)$ satisfies the compatibility equations for $\M^n(\kappa)\times\R$ if and only if for any $X,Y,Z\in\Gamma(TM)$ the two equations
\begin{align}\label{Gauss}
R(X,Y)Z=&\langle A(X),Z\rangle A(Y)-\langle A(Y),Z\rangle A(X)\\
&+\kappa\Big(\langle X,Z\rangle Y-\langle Y,Z\rangle X-\langle Y,T\rangle \langle X,Z\rangle T\nonumber\\
&-\langle X,T\rangle \langle Z,T\rangle Y+\langle X,T\rangle \langle Y,Z\rangle T
+\langle Y,T\rangle \langle Z,T\rangle X\Big),\nonumber
\end{align}
\beqt\label{Codazzi}
\nabla_XA(Y)-\nabla_YA(X)-A[X,Y]=\kappa f(\langle Y,T\rangle X-\langle X,T\rangle Y)
\eeqt
and equations \eqref{FirstequProductsp} and \eqref{SecondequProductsp} are satisfied.
\end{defn}
\begin{rem}
The relations (\ref{Gauss}) and (\ref{Codazzi}) are the Gauss and Codazzi-Mainardi equations for an isometric immersion into $\M^n(\kappa)\times\R$.
\end{rem}
Finally, we recall a result of B. Daniel (\cite{Dan}) which gives a necessary and sufficient condition for the existence of an isometric immersion of an oriented, simply connected surface $M$ into $\Ss^n(\kappa)\times\R$ or $\HH^n(\kappa)\times\R$.
\begin{thrm}[Daniel \cite{Dan}]\label{dan2}
Let $(M,\langle .,.\rangle )$ be an oriented, simply connected Riemannian manifold and $\nabla$ its Riemannian connection. Let $A$ be a field of symmetric endomorphisms $A_y:T_yM\longrightarrow T_yM$,
$T$ a vector field on $M$ and $f$ a smooth function on $M$, such that $||T||^2+f^2=1$. If $(M,\langle .,.\rangle ,A,T,f)$ satisfies the compatibility equations for $\M^n(\kappa)\times\R$, then, there exists an isometric immersion
$$F:M\longrightarrow\M^n(\kappa)\times\R$$ so that the Weingarten operator of the immersion related to the normal $\nu$ is
$$dF\circ A\circ dF^{-1}$$
and such that
$$\dt=dF(T)+f\nu.$$
Moreover, this immersion is unique up to a global isometry of
$\M^n(\kappa)\times\R$ which preserves the orientation of $\R$.
\end{thrm}
\section{Isometric immersions via spinors}
\subsection{Generalized Killing spinors}
\paragraph{The case of space forms}
We introduce the notion of generalized Killing spinors corresponding to
hypersurfaces of the space forms $\mathbb{M}^n(\kappa)$. These spinors are obtained by restriction (using (\ref{gaussspin})) of a parallel ({\it resp.} real Killing or imaginary Killing) spinor field of the ambient space $\R^n$ ({\it resp.} $\Ss^n(\kappa)$ or $\HH^n(\kappa)$). If $n$ is odd, they are the restriction of the positive part of the ambient spinor fields. We set $\eta\in\C$ such that $\kappa=4\eta^2$.
\begin{defn} A  generalized
Killing spinor on a Riemannian spin manifold $M$ with spin
connection $\nabla^{\Sigma M}$ is a solution $\varphi$ of the
 generalized Killing equation
\begin{eqnarray}\label{general_Killing_eq}\nabla^{\Sigma
M}_X\varphi=\frac{1}{2}A(X)\cdot\varphi+\eta
X\cdot \omega^{\mathbb{C}}_{n}\cdot\varphi, \end{eqnarray} for all
$X\in\Gamma(TM)$, where
 $A$ is a field of $g$-symmetric endomorphisms and
$\eta\in\mathbb{C}$. Here, $\omega^{\mathbb{C}}_{n}$ stands for the
complex volume element and "$\cdot$" is the Clifford multiplication
on $M$.
\end{defn}
\begin{rem}
Note that the complex number $\eta$ must be either real or purely imaginary because of the following well-known property of Killing spinors. If $\varphi$ satisfies
$$\nabla^{\Sigma M}_X\varphi=\eta X\cdot\varphi,$$
for all $X\in\Gamma(TM)$ then $\eta$ is either real or purely
imaginary.
\end{rem}
The norm of a generalized Killing spinor field satisfies the following
\begin{lem}\label{norm}
Let $\varphi$ be a generalized Killing spinor. Then
\begin{itemize}
\item[1.] If $\eta\in\mathbb{R}$, we have $|\varphi|=Const$.
\item[2.] If $\eta\in i\mathbb{R}$, we have
$X|\varphi|^2=-2i\eta\langle iX\cdot
\omega_n^{\C}\cdot\varphi,\varphi\rangle,$ for all $X\in\Gamma(TM)$
\end{itemize}
\end{lem}
\pf
First, we recall the well-known following lemma.
\begin{lem}\label{pdtscalaireformesreelles}
 Let $\psi$ be a spinor field and $\beta$ a real 1-form or 2-form. Then
 $$\pre\left\langle \beta\cdot\psi,\psi\right\rangle =0.$$
 \end{lem}
\noindent
Now, from this lemma, we deduce easily the proof of Lemma \ref{norm}
\begin{itemize}
\item[1.]If $\eta\in\mathbb{R}$, we have,
 $$X|\varphi|^2=2\langle\nabla_X^{\Sigma N}\varphi,\varphi\rangle=
2\langle \eta X\cdot_N\varphi,\varphi\rangle=-2
\eta\langle\varphi,X\cdot_N\varphi\rangle=0$$ and consequently
$|\varphi|=Const$.
\item[2.] If $\eta\in i\mathbb{R}$, we have
$$X|\varphi|^2=2\langle\eta
X\cdot \omega^{\mathbb{C}}_{n}\varphi,\varphi\rangle+
\langle A(X)\cdot\varphi,\varphi\rangle=-i2\eta\langle
iX\cdot \omega^{\mathbb{C}}_{n}\varphi,\varphi\rangle.$$\qed
\end{itemize}
\paragraph{The case of product spaces}

We give the following definition of the generalized Killing spinor fields corresponding to hypersurfaces of $\M^n(\kappa)\times\R$. These spinors are obtained by restriction of particular spinor fields on $\M^n(\kappa)\times\R$ playing the role of Killing spinors on space forms (see \cite{Roth4} for details). We set $\eta\in\C$ such that $\kappa=4\eta^2$.

\begin{defn}
A spinor field which satisfies the equation
\beq\label{killinggeneral} \nabla^{\Sigma M}_X\varphi&=&
-\frac{1}{2}A(X)\cdot\varphi+\eta X\cdot T\cdot\varphi+\eta
fX\cdot\varphi+\eta\left\langle X,T\right\rangle \varphi, \eeq for
all $X\in\Gamma(TM)$ where "$\cdot$'' stands for the Clifford
multiplication on $M$, $T$ is a vector field over $M$ and $f$ a
smooth function on $M$. Such a spinor field is called a generalized
Killing spinor on $\M^n(\kappa)\times\R$.
\end{defn}
\noindent
These spinor fields satisfy the following properties
\begin{prop}\label{normeconstante1}
\be[1.]
\item If $\eta\in\R$, then the norm of a generalized Killing spinor is constant.
\item If $\eta\in i\R$, then the norm of a generalized Killing spinor satisfies for any $X\in\Gamma(TM)$:
$$X|\varphi|^2= \pre\left\langle  iX\cdot T\cdot\varphi+ifX\cdot\varphi,\varphi\right\rangle. $$
\ee
\end{prop}
\noindent
{\it Proof:} We need to compute $X|\varphi|^2$ for $X\in\Gamma(TM)$. We have
$$X|\varphi|^2=2\pre\left\langle \nabla^{\Sigma M}_X\varphi,\varphi\right\rangle .$$
We replace $\nabla^{\Sigma M}_X\varphi$ by the expression given by
(\ref{killinggeneral}), and we use Lemma \ref{pdtscalaireformesreelles}
 to conclude that
 $$\pre\left\langle A(X)\cdot\varphi,\varphi\right\rangle =0,$$
 and
 $$\pre\left\langle fX\cdot\varphi,\varphi\right\rangle =0.$$
 By this lemma again, we see that
 $$\pre\left\langle X\cdot T\cdot\varphi,\varphi\right\rangle +\pre\left\langle \left\langle X,T\right\rangle\varphi,\varphi \right\rangle=0.$$
So $X|\varphi|^2=0$ and then $\varphi$ has constant norm.\\
If $\eta\in i\R$, an analogous computation yields the result. $\finpreuve$
\begin{rem}
In the case $\eta\in i\R$, the norm of $\varphi$ is not constant. Nevertheless, we can show that $\varphi$ never vanishes.
\end{rem}

\subsection{The main results}
Here, we state the main results of this paper. The first result
gives a characterization of hypersurfaces in 4-dimensional space
forms assuming the existence of two generalized Killing spinor
fields which are equivalently solutions of two Dirac equations. Part
of this result can be found in the thesis of the first author
\cite{ThLa}.
\begin{thm}\label{thmS}
Let $(M^3,g)$ be a 3-dimensional simply connected spin manifold,
 $H:M\longrightarrow\mathbb{R}$  a real valued function and $A$ a field of
 symmetric endomorphisms on $TM$.
The following statements are equivalent:
\begin{itemize}
\item[1.] The spinor fields $\varphi_{j}$, $j=1,2$, are non-vanishing solutions of the Dirac equations:
$$\left\lbrace
\begin{array}{l}
D\varphi_1=(\frac{3}{2}H+3\eta)\varphi_1,\\ \\
D\varphi_2=-(\frac{3}{2}H+3\eta)\varphi_2,
\end{array}
\right.$$
with
$\begin{cases}|\varphi_j|=Const\ \text{if}\ \eta\in \R,\\
  X|\varphi_j|^2=2\pre\left\langle \eta X\cdot\varphi_j,\varphi_j\right\rangle \quad\textrm{if}\quad\eta\in i\mathbb{R}.
\end{cases}$
\item[2.] The spinor fields $\varphi_j$, $j=1,2$, are non-trivial solutions of the generalized Killing
equations
$$\left\lbrace
\begin{array}{l}
\nabla^{\Sigma M}_X\varphi_1=\frac{1}{2}A(X)\cdot\varphi_1-\eta
X\cdot\varphi_1\\ \\
\nabla^{\Sigma M}_X\varphi_2=-\frac{1}{2}A(X)\cdot\varphi_2+\eta
X\cdot\varphi_2,
\end{array}
\right.$$
with $\frac{1}{2}\trace(A)=H$.
\end{itemize}
 Moreover both statements imply that
 \begin{itemize}
 \item[3.] there exists an isometric immersion
$F:M\hookrightarrow\mathbb{M}^4(\kappa)$ into the 4-dimensional
space form of curvature $\kappa=4\eta^2$ with mean curvature $H$ and
Weingarten tensor $dF\circ A\circ dF^{-1}$.
 \end{itemize}
 \end{thm}
\begin{rem}
 Note that in the case of $\R^4$, Assertion 3. is equivalent to Assertions 1. and 2. (see \cite{Mo})
\end{rem}

Now, we state the second result which gives a characterization of hypersurfaces into the 4-dimensional product spaces $\M^3(\kappa)\times\R$.
\begin{thm}\label{thmSR1}
 Let $(M^3,g)$ be a 3-dimensional simply connected spin manifold,
 $f,\,H:M\longrightarrow\mathbb{R}$  two real valued functions, $T$ a vector field and $A$ a field of
 symmetric endomorphisms on $TM$, such that
$$\left\lbrace
\begin{array}{l}
||T||^2+f^2=1,\\
\nabla_XT=fA(X),\\
df(X)=-\langle A(X),T\rangle.
\end{array}\right.
$$
The following statements are
 equivalent:
 \begin{itemize}
\item [1.] The spinor fields $\varphi_j$, $j=1,2,$ are non-vanishing solutions  of the generalized Dirac equations
$$\left\lbrace
\begin{array}{l}
D\varphi_1=\frac{3}{2}H\varphi_1-2\eta T\cdot\varphi_1- 3\eta f\varphi_1,\\ \\
D\varphi_2=-\frac{3}{2}H\varphi_2-2\eta T\cdot\varphi_2+ 3\eta f\varphi_2,
\end{array}
\right.$$ with constant norm if $\eta\in\R$ or satisfying
$X|\varphi|^2=\pre\big( iX\cdot
T\cdot\varphi+ifX\cdot\varphi,\varphi\big)$ if $\eta\in i\R$.
 \item[2.] The spinor fields $\varphi_j$, $j=1,2,$  are non-trivial solutions of the generalized
Killing equations
$$\left\lbrace
\begin{array}{l}
\nabla^{\Sigma M}_X\varphi_1=-\frac{1}{2}A(X)\cdot\varphi_1+\eta X\cdot T\cdot\varphi_1+\eta fX\cdot\varphi_1+\eta\left\langle X,T\right\rangle \varphi_1,\\ \\
\nabla^{\Sigma M}_X\varphi_2=\frac{1}{2}A(X)\cdot\varphi_2+\eta X\cdot T\cdot\varphi_2-\eta fX\cdot\varphi_2+\eta\left\langle X,T\right\rangle \varphi_2.
\end{array}
\right.$$
\end{itemize}
Moreover, both statements imply
\begin{itemize}
\item[3.] There exists an isometric immersion $F$ from $M$ into $\Ss^3(\kappa)\times\R$ ({\it resp.} $\HH^3(\kappa)\times\R$, with $\kappa=4\eta^2$) of mean curvature $H$ such that the Weingarten tensor related to the normal $\nu$ is given by
$$dF\circ A\circ dF^{-1}$$
and such that
$$\dt=dF(T)+f\nu.$$
\end{itemize}
 \end{thm}
\begin{rem}
As we will see in the proof (Lemma \ref{lemma3}), the condition of the existence of the two spinor fields $\varphi_1$ and $\varphi_2$ is equivalent to the existence of only one generalized Killing spinor field with $A$ a Codazzi tensor field.
\end{rem}
\section{Proof of the theorems}
We will prove Theorems \ref{thmS} and \ref{thmSR1} jointly. For this, we need three general lemmas.
\subsection{Three main lemmas}
First, we establish the following lemma which gives the Gauss
equation from a generalized Killing spinor.
\begin{lem}\label{lemma1}
Let $(M^3,g)$ be a 3-dimensional spin manifold. Assume that there
exists a non-trivial spinor field $\varphi$ solution of the
following equation \beqt \nabla^{\Sigma
M}_X\varphi=\frac{1}{2}A(X)\cdot\varphi+\eta X\cdot
T\cdot\varphi+\eta fX\cdot\varphi+\eta\left\langle X,T\right\rangle
\varphi, \eeqt where $A$, $T$ and $f$ satisfy $$\nabla_XT=fA(X),\quad df(X)=-\langle A(X),T\rangle\quad\text{and}$$
$$d^{\nabla}A(X,Y)=4\eta^2f\big(\left\langle Y,T\right\rangle X-\left\langle X,T\right\rangle Y\big),$$
then the curvature tensor $R$ of $(M,g)$ is given by
\begin{align}\label{gauss}
R(X,Y)Z=&\langle A(X),Z\rangle A(Y)-\langle A(Y),Z\rangle A(X)\\
&+\kappa\Big(\langle X,Z\rangle Y-\langle Y,Z\rangle X-\langle Y,T\rangle \langle X,Z\rangle T\nonumber\\
&-\langle X,T\rangle \langle Z,T\rangle Y+\langle X,T\rangle \langle
Y,Z\rangle T +\langle Y,T\rangle \langle Z,T\rangle X\Big).\nonumber
\end{align}
\end{lem}
\noindent {\it Proof:} We compute the spinorial curvature
$\mathcal{R}(X,Y)\varphi=\nabla^{\Sigma M}_X\nabla^{\Sigma M}_Y
\varphi- \nabla^{\Sigma M}_Y\nabla^{\Sigma M}_X
\varphi-\nabla^{\Sigma M}_{[X,Y]}\varphi$. From \cite{thRoth,Roth4},
we now that \beQ
\nabla^{\Sigma M}_X\nabla^{\Sigma M}_Y\varphi&=&\underbrace{\eta fY\cdot A(X)\cdot\varphi}_{\alpha_1(X,Y)}+\underbrace{\eta^2Y\cdot T\cdot X\cdot T\cdot\varphi}_{\alpha_2(X,Y)}+\underbrace{\eta^2fY\cdot T\cdot X\cdot\varphi}_{\alpha_3(X,Y)}\\
&&-\underbrace{\frac{\eta}{2}Y\cdot T\cdot A(X)\cdot \varphi}_{-\alpha_{4}(X,Y)}-\underbrace{\eta\left\langle A(X),T\right\rangle Y\cdot\varphi}_{-\alpha_{5}(X,Y)}+\underbrace{\eta^2fY\cdot X\cdot T\cdot\varphi}_{\alpha_{6}(X,Y)}
\eeQ
\beQ
&&+\underbrace{\eta^2\left\langle X,T\right\rangle Y\cdot T\cdot\varphi}_{\alpha_{7}(X,Y)}+\underbrace{\eta^2f^2Y\cdot X\cdot\varphi}_{\alpha_{8}(X,Y)}+\underbrace{\eta^2f\left\langle X,T\right\rangle Y\cdot\varphi}_{\alpha_{9}(X,Y)}\\
&&-\underbrace{\frac{\eta}{2}fY\cdot A(X)\cdot\varphi}_{-\alpha_{10}(X,Y)}+\underbrace{\eta f\left\langle Y,A(X)\right\rangle\varphi }_{\alpha_{11}(X,Y)}+\underbrace{\eta^2\left\langle Y,T\right\rangle X\cdot T\cdot\varphi}_{\alpha_{12}(X,Y)}
\eeQ
\beQ
&&+\underbrace{\eta^2f\left\langle Y,T\right\rangle X\cdot\varphi}_{\alpha_{13}(X,Y)}+\underbrace{\eta^2\left\langle X,T\right\rangle\left\langle Y,T\right\rangle \varphi}_{\alpha_{14}(X,Y)}-\underbrace{\frac{\eta}{2}\left\langle Y,T\right\rangle A(X)\cdot\varphi }_{-\alpha_{15}(X,Y)}\\
&&-\underbrace{\frac{1}{2}\nabla^{\Sigma M}_X(A(Y))\cdot\varphi}_{-\alpha_{16}(X,Y)}-\underbrace{\frac{\eta}{2}A(Y)\cdot X\cdot T\cdot\varphi}_{-\alpha_{17}(X,Y)}-\underbrace{\frac{\eta}{2}fA(Y)\cdot X\cdot\varphi}_{-\alpha_{18}(X,Y)}
\eeQ
\beQ
&&-\underbrace{\frac{\eta}{2}\left\langle X,T\right\rangle A(Y)\cdot\varphi}_{-\alpha_{19}(X,Y)}+\underbrace{\frac{1}{4}A(Y)\cdot A(X)\cdot\varphi}_{\alpha_{20}(X,Y)}+\underbrace{\eta\nabla^{\Sigma M}_XY\cdot T\cdot\varphi}_{\alpha_{21}(X,Y)}\\
&&+\underbrace{\eta f\nabla^{\Sigma
M}_XY\cdot\varphi}_{\alpha_{22}(X,Y)}+\underbrace{\eta\left\langle
\nabla^{\Sigma M}_XY,T\right\rangle\varphi }_{\alpha_{23}(X,Y)}.
\eeQ That is,
$$\nabla^{\Sigma M}_X\nabla^{\Sigma M}_Y\varphi=\sum_{i=1}^{23}\alpha_i(X,Y).$$
By symmetry, it is obvious that
$$\nabla^{\Sigma M}_Y\nabla^{\Sigma M}_X\varphi=\sum_{i=1}^{23}\alpha_i(Y,X).$$
On the other hand, we have
\beQ
\nabla^{\Sigma M}_{[X,Y]}\varphi&=&\underbrace{\eta[X,Y]\cdot T\cdot\varphi}_{\beta_1([X,Y])}+\underbrace{\eta f[X,Y]\cdot\varphi}_{\beta_2([X,Y])}\\
&&+\underbrace{\eta\left\langle
[X,Y],T\right\rangle\varphi}_{\beta_3([X,Y])}-\underbrace{\frac{1}{2}A[X,Y]\cdot\varphi}_{-\beta_4([X,Y])}.
\eeQ Since the connection $\nabla$ is torsion-free, we
have \beQ
\alpha_{21}(X,Y)-\alpha_{21}(Y,X)-\beta_1([X,Y])=0,\\
\alpha_{22}(X,Y)-\alpha_{22}(Y,X)-\beta_2([X,Y])=0,\\
\alpha_{23}(X,Y)-\alpha_{23}(Y,X)-\beta_3([X,Y])=0.
\eeQ
Moreover, lots of terms vanish by symmetry, namely $\alpha_1$, $\alpha_4$,
$\alpha_5$, $\alpha_{10}$, $\alpha_{11}$, $\alpha_{14}$, $\alpha_{15}$, $\alpha_{17}$, $\alpha_{18}$ and $\alpha_{19}$.\\
On the other hand, the terms $\alpha_2$, $\alpha_7$, $\alpha_8$ and $\alpha_{12}$ can be combined. Indeed, if we set
$$\alpha=\alpha_2+\alpha_7+\alpha_8+\alpha_{12},$$
then
$$\begin{array}{lll}
\alpha(X,Y)-\alpha(Y,X)
&=&\eta^2\Big[ f^2\left( Y\cdot X-X\cdot Y\right) +Y\cdot T\cdot X\cdot T -X\cdot T\cdot Y\cdot T\Big]\cdot\varphi\\\\
&=&\eta^2\Big[ f^2\left( Y\cdot X-X\cdot Y\right) +||T||^2\left( Y\cdot X-X\cdot Y\right)\Big]\cdot\varphi\\\\
&& -2\eta^2\left( \left\langle X,T\right\rangle Y\cdot T-\left\langle Y,T\right\rangle X\cdot T\right) \cdot\varphi.
\end{array}$$
If we set
$$\beta=\alpha_3+\alpha_6+\alpha_9+\alpha_{13},$$
we obtain
$$
\beta(X,Y)-\beta(Y,X)=\eta^2f\left( \left\langle Y,T \right\rangle
X-\Big\langle X,T \right\rangle Y\Big) \cdot\varphi.
$$
Finally, we get
\beQ
\mathcal{R}(X,Y)\varphi&=&\frac{1}{4}\left( A(Y)\cdot A(X)-A(X)\cdot A(Y)\right) \cdot\varphi-\frac{1}{2}d^{\nabla}A(X,Y)\cdot\varphi\\\\
&&+\eta^2f\left( \left\langle Y,T \right\rangle X-\left\langle X,T \right\rangle Y\right) \cdot\varphi+\eta^2\left( Y\cdot X-X\cdot Y\right)\cdot\varphi \\\\
&&-2\eta^2\left( \left\langle X,T\right\rangle Y\cdot T-\left\langle
Y,T\right\rangle X\cdot T\right) \cdot\varphi. \eeQ Since we assume
that $A$ satisfies the following Codazzi equation
$$d^{\nabla}A(X,Y)=4\eta^2f\big(\left\langle Y,T\right\rangle X-\left\langle X,T\right\rangle Y\big),$$
we have
\beq\label{integrability_condition}
~~~~~~\mathcal{R}(X,Y)\varphi&=&\frac{1}{4}\left( A(Y)\cdot A(X)-A(X)\cdot A(Y)\right)\cdot\varphi\\
&&+\eta^2f\left( \left\langle Y,T \right\rangle X-\left\langle X,T
\right\rangle Y\right) \cdot\varphi+\eta^2\left( Y\cdot X-X\cdot
Y\right)\cdot\varphi \nonumber\eeq Now, let $X=e_i$ and $Y=e_j$ with
$i\neq j$. The Ricci identity sA(Y)s that:
\begin{eqnarray}\label{Ricciidentity}
\mathcal{R}(e_i, e_j)\cdot\varphi=\frac{1}{2}[R_{ijik} e_j-R_{ijij}
e_k-R_{ijjk} e_i]\cdot\varphi,
\end{eqnarray}
where $(i, j, k)$ is any cyclic permutation of $(1, 2, 3)$. \\\\
Further with a simple computation we find
\begin{eqnarray*}
A(e_j)\cdot A(e_i)-A(e_i)\cdot A(e_j)&=&2(A_{ik}A_{jj}-A_{ij}A_{jk})e_i\\
&&-2(A_{ik}A_{ji}-A_{ii}A_{jk})e_j+2(A_{ij}A_{ji}-A_{ii}A_{jk})e_k.
\end{eqnarray*}
With the integrability condition (\ref{integrability_condition})
this yields
\begin{eqnarray*}
(\nabla_{e_j}A)(e_i)-(\nabla_{e_i}A)(e_j)&=&(R_{ijjk}-(A_{ik}A_{jj}-A_{ij}A_{jk})+\kappa f^2)e_i\\
&&-(R_{ijik}-(A_{ik}A_{ji}-A_{ii}A_{jk})+\kappa f^2)e_j\\
&&
+(R_{ijij}-(A_{ij}A_{ji}-A_{ii}A_{jk})+\kappa f^2)e_k\\
&&+\kappa f\left( \left\langle e_i,T\right\rangle e_j-\left\langle e_i,T\right\rangle e_i \right),
\end{eqnarray*}
which proves that, if $A$ is a Codazzi tensor, it satisfies the
Gauss equation too. This observation was made by Morel (\cite{Mo})
in the Riemannian case for a parallel tensor $A$. We point out that
the converse is also true.
$\finpreuve$\\
Now, we state a second lemma which will give the equivalence between
the Dirac equation and the Killing equation (up to a condition on
the norm of the spinor field).
\begin{lem}\label{lemma2}
Let $(M^3,g)$ be a 3-dimensional spin manifold. Assume that there
exists a non-trivial spinor field $\varphi$, solution of the
following equation \beqt\label{eqdirac}
D\varphi=\frac{3}{2}H\varphi-2\eta T\cdot\varphi-3\eta f\varphi,
\eeqt where the norm of $\varphi$ satisfies for all $X\in\Gamma(TM)$
$$X|\varphi|^2=2\pre\left\langle \eta X\cdot T\cdot\varphi+\eta fX\cdot\varphi,\varphi\right\rangle .$$
Then $\varphi$ is a solution of the following generalized Killing
spinors equation
 \beqt \nabla^{\Sigma
M}_X\varphi=\frac{1}{2}A(X)\cdot\varphi+\eta X\cdot T\cdot\varphi+\eta
fX\cdot\varphi+\eta\left\langle X,T\right\rangle \varphi. \eeqt
\end{lem}
\noindent
{\it Proof:}
The 3-dimensional complex spinor space is
$\Sigma_3\cong\mathbb{C}^{2}$. The complex spin representation is
then real 4-dimensional. We now define the map
\function
{f}{\mathbb{R}^3\oplus\mathbb{R}}{\Sigma_3}{(v,~r)}{v\cdot\varphi+ r\varphi,}
where $\varphi$ is a given non-vanishing spinor.\\
Obviously $f$ is an isomorphism. Then for all $\psi\in\Sigma_3$
there is a unique pair
$(v,~r)\in(\mathbb{R}^3\oplus\mathbb{R})\cong
T_pM^3\oplus\mathbb{R}$, such that $\psi=v\cdot\varphi+r\varphi$.\\
Consequently $\big({\nabla^{\Sigma M}_X\varphi}\big)_p\in
\Gamma(T_p^*M\otimes\Sigma_3)$ can be expressed as follows:
\begin{eqnarray*}
\nabla^{\Sigma M}_X\varphi=B(X)\cdot\varphi+\omega(X)\varphi,
\end{eqnarray*}
for all $p\in M$ and for all vector fields $X$, with $\omega$
a 1-form and $B$ a (1,1)-tensor field.\\
Moreover we have
\begin{eqnarray*}
X|\varphi|^2=2\pre\langle\nabla^{\Sigma M}_X\varphi,~\varphi\rangle=2\langle \omega(X)\varphi,~\varphi\rangle\Rightarrow \omega(X)=\frac{d(|\varphi|^2)}{2|\varphi|^2}(X).
\end{eqnarray*}
which yields $\omega(X)=\pre\left\langle \eta X\cdot T\cdot\varphi+\eta fX\cdot\varphi,\dfrac{\varphi}{|\varphi|^2}\right\rangle$.\\
Now, let $B=S+U$ with $S$ the symmetric and $U$ the skew-symmetric
part of $B$. Let $\{e_i\}$ be an orthonormal basis of $TM$ and
$\varphi$ a solution of the Dirac equation (\ref{eqdirac}). We have
\begin{eqnarray*}
D\varphi&=& \sum_{i=1}^3 e_i\cdot\nabla^{\Sigma M}_{e_i}\varphi=\sum_{i,j=1}^3 e_i \cdot B_{ij}e_j\cdot\varphi+\sum_{j=1}^3\omega(e_j)e_j\cdot\varphi\\
        &=& \sum_{i=1}^3 U_{ij} e_i\cdot e_j\cdot\varphi +\sum_{i=1}^3 S_{ii} e_i\cdot e_i\cdot\varphi +
\sum_{i\neq j}^3\underbrace{S_{ij}}_{sym.} \underbrace{e_i\cdot e_j}_{skew-sym.}\cdot\varphi+W\cdot\varphi,
\end{eqnarray*}
where $W$ is the vector field defined by $W:=\sum_{j=1}^3\omega(e_j)e_j$. Then,
\begin{eqnarray*}
D\varphi&=& -2\sum_{i<j}^3 U_{ij} e_i\cdot e_j\cdot\varphi + \sum_{i=1}^3 S_{ii} e_i\cdot e_i\cdot\varphi+W\cdot\varphi\\
        &=& -2(U_{12} e_1\cdot e_2+U_{13} e_1\cdot e_3+U_{23} e_2\cdot e_3)\cdot\varphi-\trace(B)\varphi+W\cdot\varphi\\
\end{eqnarray*}

We recall that the complex volume element $\omega_3^{\C}=-e_1\cdot
e_2\cdot e_3$ acts as the identity on $\Sigma\,M$, where
$\{e_1,e_2,e_3\}$ is a local orthonormal frame of $TM$. So we deduce
that for any spinor field on $M$, $e_i\cdot e_j\cdot\varphi=e_k\cdot\varphi$, where $(i,j,k)$ is a cyclic permutation of $(1,2,3)$. From this fact, we get
\begin{eqnarray*}
D\varphi&=& -2(U_{12} e_3-U_{13} e_2+U_{23} e_1)\cdot\varphi-\trace(B)\varphi+W\cdot\varphi.\\
\end{eqnarray*}
On the other hand, we have
\begin{eqnarray*}
D\varphi&=&\frac{3}{2}H\varphi-2\eta T\cdot\varphi-3\eta f\varphi.
\end{eqnarray*}
Note that $\pre\langle (U_{12} e_3-U_{13} e_2+U_{23} e_1)\varphi,\varphi\rangle=0$ and $\pre\left\langle W\cdot\varphi,\varphi\right\rangle=0 $.
It follows that
$$
\frac{3}{2}H|\varphi|^2-2\pre\left\langle \eta T\cdot\varphi,\varphi\right\rangle -3\pre\left\langle \eta f\varphi,\varphi\right\rangle =-\trace(B)|\varphi|^2.$$
Moreover, since $\left\lbrace \frac{\varphi}{|\varphi|},\frac{e_1\cdot\varphi}{|\varphi|},\frac{e_2\cdot\varphi}{|\varphi|},\frac{e_3\cdot\varphi}{|\varphi|}\right\rbrace $ is
an orthonormal frame of $\Sigma_pM$ for the real scalar product $\left\langle .,. \right\rangle $, we deduce that
\beQ
-2\big(U_{12} e_3-U_{13} e_2+U_{23} e_1\big)\cdot\varphi&=&-3\eta f\varphi-W\cdot\varphi-2\eta T\cdot\varphi+2\pre\left\langle \eta T\cdot\varphi,\varphi\right\rangle \varphi\\
&& +3\pre\left\langle \eta f\varphi,\varphi\right\rangle\varphi.
\eeQ Further we compute
\begin{eqnarray*}
\langle U(e_j)\cdot\varphi,~e_i\cdot\varphi\rangle &=&\sum_k^3 U_{kj}
\underbrace{\langle e_k\cdot\varphi,~e_i\cdot\varphi\rangle}_{=-\langle e_i\cdot
e_k\cdot\varphi,~\varphi\rangle=0,~i\neq k}= U_{ij}|\varphi|^2.
\end{eqnarray*}
Consequently, for $i,j\in\{1,2,3\}$, we have
\beQ
-2\sum_{k<l}^3U_{lk}\left\langle e_k\cdot e_l\cdot\varphi,e_i\cdot e_j\cdot\varphi\right\rangle &=& -3\left\langle \eta f\varphi,e_i\cdot e_j\cdot\varphi \right\rangle - \left\langle W\cdot\varphi,e_i\cdot e_j\cdot\varphi\right\rangle\\
&& -2\left\langle \eta T\cdot\varphi,e_i\cdot e_j\cdot\varphi\right\rangle
+2\left\langle \eta T\cdot\varphi,\varphi\right\rangle\left\langle \varphi,e_i\cdot e_j\cdot\varphi\right\rangle \\
&&+3\left\langle\eta f\varphi,\varphi\right\rangle\left\langle
\varphi,e_i\cdot e_j\cdot\varphi\right\rangle. \eeQ Moreover, in the
3-dimensional case at most three of the four indices could be distinct.
Then, for $m\neq n$, $\langle e_m\cdot
e_n\cdot\varphi,~\varphi\rangle=0$ holds and
 as the trace of  a skew-symmetric tensor vanishes, we have:
$\langle e_k\cdot e_l\cdot\varphi,~e_j\cdot e_i\cdot\varphi\rangle\neq 0
\Leftrightarrow k=i, l=j$ or $k=j, l=i$, $i\neq j$,  which yield
\begin{eqnarray*}
-2U_{ij}|\varphi|^2&=&-2\langle U(e_j)\cdot\varphi,e_i\cdot\varphi\rangle\\\\
&=&-3\left\langle \eta f\varphi,e_i\cdot e_j\cdot\varphi \right\rangle - \left\langle W\cdot\varphi,e_i\cdot e_j\cdot\varphi\right\rangle -2\left\langle \eta T\cdot\varphi,e_i\cdot e_j\cdot\varphi\right\rangle\\\\
&&+3\left\langle \eta f\varphi,\varphi\right\rangle \left\langle e_j\cdot\varphi,e_i\cdot\varphi\right\rangle +2\left\langle\eta T\cdot\varphi,\varphi \right\rangle \left\langle e_j\cdot\varphi,e_i\cdot\varphi\right\rangle.
\end{eqnarray*}
Then, we deduce that
\beq\label{U}
-2U(X)&=&X\cdot W\cdot\varphi-\left\langle X\cdot W\cdot\varphi,\varphi\right\rangle\dfrac{\varphi}{|\varphi|^2}-2\eta X\cdot T\cdot\varphi\nonumber\\
&&+2\left\langle \eta X\cdot T\cdot\varphi,\varphi\right\rangle \dfrac{\varphi}{|\varphi|^2}+3\left\langle \eta f\varphi,\dfrac{\varphi}{|\varphi|^2}\right\rangle X\cdot\varphi\nonumber\\
&&+2\left\langle \eta
T\cdot\varphi,\dfrac{\varphi}{|\varphi|^2}\right\rangle
X\cdot\varphi-3\eta fX\cdot\varphi+3\left\langle \eta
fX\cdot\varphi,\varphi \right\rangle \dfrac{\varphi}{|\varphi|^2}.
\eeq From now on, we will consider separately the cases
$\eta\in\R$ and $\eta\in i\R$.
\subsubsection*{The case $\eta\in\R$}
Since $\eta$ is real, the norm of $\varphi$ is constant and so $\omega(X)=0$ for any vector field $X$. Consequently, using Lemma \ref{pdtscalaireformesreelles}, we get
\beQ
U(X)\cdot\varphi&=& \eta X\cdot T\cdot\varphi-\eta\left\langle X\cdot T\cdot\varphi,\varphi\right\rangle \dfrac{\varphi}{|\varphi|^2}\\
&=& \eta X\cdot T\cdot\varphi +\eta\left\langle
X,T\right\rangle \varphi. \eeQ Moreover,
\begin{eqnarray*}
Q_{\varphi}(e_i,e_j)&=&\frac{1}{2}\left\langle e_i\cdot\nabla^{\Sigma M}_{e_j}\varphi+e_j\cdot\nabla^{\Sigma M}_{e_i}\varphi,\frac{\varphi}{|\varphi|^2}\right\rangle\\
                    &=&\frac{1}{2}\left\langle\sum_{k}^3 S_{jk} e_i \cdot e_k\cdot\varphi+\sum_{k}^3
S_{ik}e_j\cdot e_k\cdot\varphi, \frac{\varphi}{|\varphi|^2}\right\rangle
\\
&=&-S_{ij}|\varphi|^2 ~~~\Rightarrow
S(X)=-Q_{\varphi}(X).
\end{eqnarray*}
Now, we set
$$A(X):=2Q_{\varphi}(X)+2\eta fX.$$
Finally, we obtain \beqt \nabla^{\Sigma
M}_X\varphi=\frac{1}{2}A(X)\cdot\varphi+\eta X\cdot
T\cdot\varphi+\eta fX\cdot\varphi+\left\langle X,T\right\rangle
\varphi, \eeqt which achieves the proof in the case $\eta\in\R$.
\subsubsection*{The case $\eta\in i\R$}

Here, $\eta$ is not real and so the norm of $\varphi$ is not constant but satisfies
$$X|\varphi|^2=2\pre\left\langle \eta X\cdot T\cdot\varphi+\eta fX\cdot\varphi,\varphi\right\rangle .$$
Then
\beqt\label{omega}
\omega(X)=\dfrac{X|\varphi|^2}{2|\varphi|^2}=\dfrac{1}{2|\varphi|^2}\pre\left\langle \eta X\cdot T\cdot\varphi+\eta fX\cdot\varphi,\varphi\right\rangle.
\eeqt
Like in the case $\eta\in\R$, we have $S(X)=-Q_{\varphi}(X)$ and we set
$$A(X):=2Q_{\varphi}(X)+V(X),$$
where $V(X)$ is the symmetric endomorphism field defined by
\beq\label{V}
V(X,Y)&=&2\pre\left\langle \eta\left\langle X,Y\right\rangle T\cdot\varphi,\varphi\right\rangle +2\pre\left\langle \eta f\left\langle X,Y\right\rangle\varphi,\varphi\right\rangle\nonumber\\
&&+\pre\left\langle  \eta\left( \left\langle X,T\right\rangle Y+
\left\langle Y,T\right\rangle X\right)\cdot\varphi,\varphi
\right\rangle. \eeq Since $$\nabla^{\Sigma
M}_X\varphi=S(X)\cdot\varphi+U(X)\cdot\varphi+\omega(X)\varphi,$$ we
deduce from (\ref{omega}), (\ref{U}) and (\ref{V}) that \beqt
\nabla^{\Sigma M}_X\varphi=\frac{1}{2}A(X)\cdot\varphi+\eta X\cdot
T\cdot\varphi+\eta fX\cdot\varphi+\eta\left\langle
X,T\right\rangle \varphi. \eeqt $\finpreuve$

Now, we give a final lemma which will allow us to use Lemma
\ref{lemma1} for the proof of Theorems \ref{thmS} and \ref{thmSR1}.
Indeed, in Theorems \ref{thmS} and \ref{thmSR1}, we do not suppose
anything about the symmetric tensor $A$. Nevertheless, the existence
of two generalized Killing spinor fields implies that $A$ is
Codazzi.

\begin{lem}\label{lemma3}
Let $(M^3,g)$ a 3-dimensional spin manifold. Assume that there exist
two non-trivial spinor fields $\varphi_1$ and $\varphi_2$ such that \beqt
\nabla^{\Sigma}_X\varphi_1=\frac{1}{2}A(X)\cdot\varphi_1+\eta X\cdot
T\cdot\varphi_1+\eta fX\cdot\varphi_1+\left\langle X,T\right\rangle
\varphi_1, \eeqt and \beqt \nabla^{\Sigma
M}_X\varphi_2=-\frac{1}{2}A(X)\cdot\varphi_2+\eta X\cdot T\cdot\varphi_2-\eta
fX\cdot\varphi_2+\left\langle X,T\right\rangle \varphi_2, \eeqt where $A$, $T$
and $f$ satisfy $$\nabla^{\Sigma M}_XT=fA(X),\quad df(X)=-\langle
A(X),T\rangle,$$ then the tensor $A$ satisfies the Codazzi-Mainardi
equations, that is
$$d^{\nabla}A(X,Y)=4\eta^2f\big(\left\langle Y,T\right\rangle X-\left\langle X,T\right\rangle Y\big).$$
\end{lem}
\pf
From the proof of Lemma \ref{lemma1}, we know that the equation satisfied by $\varphi_1$ implies
\begin{eqnarray}\label{CG1}
(\nabla_{e_j}A)(e_i)-(\nabla_{e_i}A)(e_j)&=&(R_{ijjk}-(A_{ik}A_{jj}-A_{ij}A_{jk})+\kappa f^2)e_i\nonumber\\
&&-(R_{ijik}-(A_{ik}A_{ji}-A_{ii}A_{jk})+\kappa f^2)e_j\nonumber \\
&&
+(R_{ijij}-(A_{ij}A_{ji}-A_{ii}A_{jk})+\kappa f^2)e_k\\
&&+\kappa f\left( \left\langle e_i,T\right\rangle e_j-\left\langle e_i,T\right\rangle e_i \right) .\nonumber
\end{eqnarray}
On the other hand, by an analogous computation for the spinor field
$\varphi_2$, we get
\begin{eqnarray*}
-(\nabla_{e_j}A)(e_i)+(\nabla_{e_i}A)(e_j)&=&(R_{ijjk}-(A_{ik}A_{jj}-A_{ij}A_{jk})+\kappa f^2)e_i\\
&&-(R_{ijik}-(A_{ik}A_{ji}-A_{ii}A_{jk})+\kappa f^2)e_j\\
&&
+(R_{ijij}-(A_{ij}A_{ji}-A_{ii}A_{jk})+\kappa f^2)e_k\\
&&-\kappa f\left( \left\langle e_i,T\right\rangle e_j-\left\langle
e_i,T\right\rangle e_i \right).
\end{eqnarray*}
If we combine the last two equalities, we get
$$\left\lbrace
\begin{array}{l}
R_{ijjk}-(A_{ik}A_{jj}-A_{ij}A_{jk})+\kappa f^2=0,\\
R_{ijik}-(A_{ik}A_{ji}-A_{ii}A_{jk})+\kappa f^2=0,\\
R_{ijij}-(A_{ij}A_{ji}-A_{ii}A_{jk})+\kappa f^2=0,
\end{array}\right. $$
that is exactly the Gauss equation. Then, we get immediately from equation \eqref{CG1} that $A$ also satisfies the Codazzi equation
$$d^{\nabla}A(X,Y)=4\eta^2f\big(\left\langle Y,T\right\rangle X-\left\langle X,T\right\rangle Y\big),$$
for all vector fields $X$ and $Y$.$\finpreuve$

\subsection{Proof of the Theorems}\label{proof}
The proof of the theorems follows easily from Lemmas \ref{lemma1}, \ref{lemma2} and \ref{lemma3} with
$$\left\lbrace
\begin{array}{ll}
\eta=0&\text{for}\;\R^4,\\\\
\eta=\dfrac{1}{2},\;T=0,\;f=1&\text{for}\;\Ss^4,\\\\
\eta=\dfrac{i}{2},\;T=0,\;f=1&\text{for}\;\HH^4,\\\\
\eta=\dfrac{1}{2}&\text{for}\;\Ss^3\times\R,\\\\
\eta=\dfrac{i}{2}&\text{for}\;\HH^3\times\R.
\end{array}
\right. $$

 Indeed, Lemma \ref{lemma2} gives the equivalence between Assertions $1.$ and $2.$ of the theorems,
 that is, between the existence of a generalized Killing spinor  and a Dirac spinor satisfying an additional norm condition.\\
\indent The proof of $2.\Longrightarrow 3.$ is an immediate
consequence of Lemmas \ref{lemma1} and \ref{lemma3}. From Lemma
\ref{lemma3}, the problem is reduced to the case of only one
generalized Killing spinor field, but with $A$ a Codazzi tensor.
Now, if the tensor $A$ satisfies the Codazzi-Mainardi equation, then
by Lemma \ref{lemma1}, it satisfies also the Gauss equation. It is
well-known that if the Gauss and Codazzi-Mainardi equations are
satisfied for a simply connected manifold, then it can be immersed
isometrically in the corresponding space form. For the case of
product spaces, by the result of Daniel (\cite{Dan}), to get an
isometric immersion, the two additional conditions (\ref{FirstequProductsp}) and
(\ref{SecondequProductsp}) are needed.$\finpreuve$

\begin{rem}
Conversely, the existence of one generalized Killing spinor field
$\varphi_1$ with Codazzi tensor field $A$ implies the existence of a
second spinor field $\varphi_2$. Indeed, as we just saw, $M$ is
isometrically immersed into $\M^4(\kappa)$ or
$\M^3(\kappa)\times\R$. Then, one just defines $\varphi_2$ as
$\nu\cdot\varphi_1$, where $\nu$ is the normal unit vector field.
Thus, if $\varphi_1$ satisfies
$$\nabla^{\Sigma M}_X\varphi_1=-\frac{1}{2}A(X)\cdot\varphi_1+\eta X\cdot T\cdot\varphi_1+\eta fX\cdot\varphi_1+\eta\left\langle X,T\right\rangle \varphi_1,$$
then, by a straightforward computation, $\varphi_2$ satisfies
$$\nabla^{\Sigma M}_X\varphi_2=\frac{1}{2}A(X)\cdot\varphi_2+\eta X\cdot T\cdot\varphi_2-\eta fX\cdot\varphi_2+\eta\left\langle X,T\right\rangle \varphi_2.$$

\end{rem}

\section{Application: Non-existence of isometric immersions for 3-dimensional geometries}
In  \cite{Pet} and \cite{Mas}, for instance,  it is shown that there exist no isometric immersions for
certain 3-dimensional homogeneous spaces into the Euclidean 4-space.
As an application of Theorem \ref{thmS} we give a short
non-spinorial proof of the non-existence of such immersions for
certain three-dimensional $\eta-$manifolds including the above
homogeneous spaces.
\subsection{Preliminaries the some 3-dimensional geometries}
In this section, we will give some basic facts about 3-dimensional
homogeneous manifolds. A complete description can be found in
\cite{Sco}. Let $(M^3,g)$ be a 3-dimensional Riemannian homogeneous
manifold. We denote by $d$ the dimension of its isometry group. The
possible values of $d$ are $3$, $4$ and $6$. If $d$ is equal to $6$,
then $M$ is a space form $\mathbb{M}^3(\kappa)$. There is only one
geometry with $d$ equal to $3$, namely, the solvable group $Sol_3$.
Finally, if $d=4$, then, there are $5$ possible models.
\subsubsection{The manifolds $\Ekt$ with $\tau\neq0$}
Such manifolds are Riemannian fibrations over $2$-dimensional space forms. They are denoted by $\Ekt$ where $\kappa$ is the curvature of the base of the fibration and $\tau$ is the bundle curvature, that is the defect for the fibration to be a product. Note that $\kappa\neq4\tau^2$, if not, the manifold is a space form.  Table 1. gives the classification of these possible geometries.
\begin{center}
\begin{tabular}{|c|c|c|c|}
\hline
&$\kappa>0$&$\kappa=0$&$\kappa<0$\\
\hline
$\tau=0$&$\Ss^2(\kappa)\times\R$& $\R^3$ &$\HH^2(\kappa)\times\R$\\
\hline
$\tau\neq0$&$(\Ss^3,g_{Berger})$&$Nil_3$&$\widetilde{PSL_2(\R)}$\\
\hline
\end{tabular}\\
\vspace{0.2cm}
Table 1: \textit{Classification of }$\Ekt$
\end{center}
\vspace{0.3cm}
From now on, we will focus on the non-product case, {\it i.e.}, $\tau\neq0$. In this case, $\Ekt$ carries a unitary Killing vector field $\xi$ tangent to the fibers and satisfying
$\nabla_X\xi=\tau X\wedge\xi.$
Moreover, there exits a direct local orthonormal frame $\{e_1,e_2,e_3\}$ with $e_3=\xi$ and  such that the Christoffel symbols are
\begin{equation}\label{christoffel}
\left\lbrace
\begin{array}{l}
\Gamma_{12}^3=\Gamma_{23}^1=-\Gamma_{21}^3=-\Gamma_{13}^2=\tau,\\
\Gamma_{32}^1=-\Gamma_{31}^2=\tau-\dfrac{\kappa}{2\tau}, \\
\Gamma_{ii}^i=\Gamma_{ij}^i=\Gamma_{ji}^i=\Gamma_{ii}^j=0,\quad\forall\,i,j\in\{1,2,3\}.
\end{array}
\right.
\end{equation}
In particular, we deduce from these Christoffel symbols that $\Ekt$ is $\eta$-Einstein. Precisely, we have
$$Ric=
\left(
\begin{array}{ccc}
 \kappa-2\tau^2 & 0  & 0  \\
 0 &  \kappa-2\tau^2 & 0  \\
 0 &  0 & 2\tau^2
\end{array}
\right)
$$
in the local frame $\{e_1,e_2,\xi\}$. Moreover, from this and the local expression of the spinorial Levi-Civita connection, we deduce that there exists on $\Ekt$ a spinor field $\varphi$ satisfying
\begin{equation}\label{presquekillingekt}
 \nabla_{e_1}\varphi=\dfrac{1}{2}\tau e_1\cdot\varphi,\
 \nabla_{e_2}\varphi=\dfrac{1}{2}\tau e_2\cdot\varphi,\
\nabla_{\xi}\varphi=\dfrac{1}{2}\left(\dfrac{\kappa}{2\tau}-\tau\right)\xi\cdot\varphi.
\end{equation}
One can refer to \cite{Roth4} for details.
\subsubsection{The Lie group $Sol_3$}
The solvable Lie group $Sol_3$ is the semi-direct product $\R^2\rtimes\R$ where $t\in\R$ acts on $\R^2$ by the transformation $(x,y)\longrightarrow(e^tx,e^ty)$. Then, we identify $Sol_3$ with $\R^3$ and the group multiplication is defined by
$$ (x,y,z)\cdot(x',y',z')=(x+e^{-z}x',y+e^{z}y',z+z').$$
The frame $e_1=e^{-z}\partial_x$, $e_2=e^z\partial_y$, $e_3=\partial_z$ is orthonormal for the left-invariant metric
$$ds^2=e^{2z}dx^2+e^{-2z}dy^2+dz^2.$$
We easily check that in the frame $\{e_1,e_2,e_3\}$, the Christoffel symbols are
$$\Gamma_{11}^3=\Gamma_{23}^2=-\Gamma_{13}^1=-\Gamma_{22}^3=-1,$$
and the other identically vanish. So, we deduce the existence of a special spinor field $\varphi$ on $Sol_3$ satisfying
  \begin{equation}\label{presquekillingsol}
 \overline{\nabla}_{e_1}\varphi=\dfrac{1}{2} e_2\cdot\varphi,\
 \overline{\nabla}_{e_2}\varphi=\dfrac{1}{2} e_1\cdot\varphi,\
\overline{\nabla}_{\xi}\varphi=0,
\end{equation}
and the Ricci curvature in the frame $\{e_1,e_2,e_3\}$ is given by
$$
\left(
\begin{array}{ccc}
 0 & 0  & 0  \\
 0 &  0 &  0 \\
 0 & 0  & -2
\end{array}
\right).
$$
Details can be found in \cite{Hab}.
\subsubsection{The hyperbolic fibration $\mathbb{T}_B^3$}
This last example is the hyperbolic fibration defined in \cite{Mey}. Let
$B$ be a matrix of $\mathrm{SL}_2(\mathbb{Z})$, which can be
considered as a diffeomorphism of the flat torus $\mathbb{T}^2$ and
admit two eigenvalues $\alpha$ and $\frac{1}{\alpha}$. Now let
$\mathbb{T}_B^3$ be the 3-dimensional manifold defined by
$\mathbb{T}_B^3=\mathbb{T}^2\times\mathbb{R} / \equiv$, where
$\equiv$ is the equivalence relation defined by
$(x,y)\equiv(B(x),y+1)$. We denote by $b$ the slope of the
eigenvector associated to the eigenvalue $\frac{1}{\alpha}$. Thus,
$\mathbb{T}_B^3$ is a compact manifold of universal covering
$\mathbb{R}^3$ equipped with a Riemannian metric for which the base
$\{e_1,e_2,e_3\}$ defined as follows is orthonormal
$$e_1=\alpha^{-z}(-b\partial_x+\partial_y),\ e_2=\alpha^{z}(\partial_x+b\partial_y),\ e_3=\partial_z.$$
One can easily check that
$$[e_1,e_2]=0,\ [e_1,e_3]=\ln(\alpha)e_1,\ [e_2,e_3]=-\ln(\alpha)e_2,$$
and that the Christoffel symbols are given by
$$\Gamma_{11}^3=\Gamma_{23}^2=-\Gamma_{13}^1=-\Gamma_{22}^3=-\ln(\alpha),$$
with the other identically zero. The Ricci curvature is given by the following matrix in the frame $\{e_1,e_2,e_3\}$
$$
\left(
\begin{array}{ccc}
 0 & 0  & 0  \\
 0 &  0 &  0 \\
 0 & 0  & -2\ln(\alpha)^2
\end{array}
\right).
$$
From the expression of the Christoffel symbols, there exists a spinor field $\varphi$ satisfying
\beqt\label{spinorTB3}
\nabla_{e_1}\varphi=\frac{1}{2}\ln(\alpha)e_2\cdot\varphi,\ \nabla_{e_2}\varphi=\frac{1}{2}\ln(\alpha)e_1\cdot\varphi,\ \nabla_{e_3}\varphi=0.\eeqt
\subsection{A non-existence result}
Here is the main result of this section.
 \begin{prop} \label{Prop_3homspaces}The 3-dimensional manifolds $Nil_3$,
$Sol_3$, $\widetilde{PSl_2(\mathbb{R})}$, the Berger spheres and the tori $\mathbb{T}_B^3$ cannot be isometrically
immersed into $\mathbb{R}^4$, even locally.
\end{prop} We start by giving the following
 \begin{lem} \label{Linalg_Lemma} Let $(M^3,g)$ be an oriented Riemannian manifold which is $\eta-$Einstein, i.e. $Ric=\lambda g + \eta\xi\otimes\xi$, with $\eta\neq0$. Assume that there exists a
non-trivial spinor field $\varphi$ such that
$\nabla^{\Sigma
M}_X\varphi=-\frac{1}{2}A(X)\cdot\varphi$, where $A$ is a symmetric
endomorphism field. Then,
\begin{itemize}
\item[1.] If $\lambda\neq-\eta$, and A is Codazzi, then
$$A=\pm\left(
      \begin{array}{ccc}
        \sqrt{\frac{|\lambda+\eta|}{2}} & 0 & 0 \\
        0 & \sqrt{\frac{|\lambda+\eta|}{2}} & 0 \\
        0 & 0 & \frac{\lambda-\eta}{\sqrt{2|\lambda+\eta|}} \\
      \end{array}
    \right)
$$
in an orthonormal frame $\{e_1, e_2, \xi\}$.
\item[2.]If $\lambda=-\eta$, then $A$ cannot be Codazzi.
\item[3.] If $\lambda=0$ and $\eta<0$, then $A$ cannot be Codazzi.
\end{itemize}
 \end{lem}
 \pf Using the fact that $A$ is Codazzi, a simple calculation shows
 $$R^{\Sigma
M }(X,Y)\cdot\varphi=\frac{1}{4}(A(Y)\cdot A(X)-A(X)\cdot
A(Y))\cdot\varphi.$$ Then the Ricci identity (\ref{Ricciidentity}) yields
$$Ric(X)\cdot\varphi=\trace(A)A(X)\cdot\varphi-A^2(X)\cdot\varphi.$$
Now if the manifold is $\eta-$Einstein, we get
$$\Big(\lambda X+\eta\langle X,\xi\rangle\xi-\trace (A)A(X)+A^2(X)\Big)\cdot\varphi=0.$$
Since $\varphi$ is a non-trivial generalized Killing spinor, it
never vanishes. Consequently
\beqt
\lambda X+\eta\langle
X,\xi\rangle\xi-\trace(A)A(X)+A^2(X)=0.\label{equation_from_Ricci}
\eeqt
Let $\{e_1,e_2,e_3\}$ be a diagonalizing frame of $A$, then from
equation \eqref{equation_from_Ricci} $e_3$ can always be chosen to be
$\xi$ and $e_1$, $e_2$ orthogonal to $\xi$. Now denote by $a_1$,
$a_2$, $a_3$ the respective eigenvalues. Then equation
\eqref{equation_from_Ricci} leads to
$$\left\{
    \begin{array}{l}
      a_1a_2=\frac{\lambda+\eta}{2}, \\
      a_2a_3=\frac{\lambda-\eta}{2},\\
      a_1a_3=\frac{\lambda-\eta}{2}. \\
    \end{array}
  \right.
$$
If $\lambda=-\eta$, then this system has no solutions. If $\lambda=0$ and $\eta<0$, then we have $a_1=a_2$, and so $a_1^2=\frac{\eta}{2}<0$, which is not possible because $a_1$ is a real number. Thus, in these two cases, $A$ cannot be Codazzi. If $\lambda\neq-\eta$ simple computations yield
the result.\qed\\\\
{\it Proof of Proposition \ref{Prop_3homspaces}:} Let $M=Nil_3$,
$Sol_3$, $\widetilde{PSl_2(\mathbb{R})}$, $\mathbb{T}_B^3$ or a
Berger sphere and assume that $M$ is isometrically immersed in
$\mathbb{R}^4$. Then there exists a spinor $\varphi$ on $M$
verifying $\nabla^{\Sigma M}_X\varphi=-\frac{1}{2}A(X)\cdot\varphi,$
where $A$ is shape operator of the immersion and hence Codazzi.
Moreover, all these manifolds are $\eta$-Einstein. For $Sol_3$ and
$\mathbb{T}_B^3$, we have $\lambda=0$ and $\eta<0$, so from Lemma
\ref{Linalg_Lemma}, $A$ cannot be Codazzi and such a spinor cannot
exist. This leads to a contradiction. In the case of $Nil_3$,
$\widetilde{PSl_2(\mathbb{R})}$ and Berger spheres, we have
$\lambda=\kappa-2\tau^2$ and $\eta=2\tau^2$. Since
$\kappa\neq4\tau^2$, then $\lambda\neq-\eta$ and $A$ is as in part 1
of Lemma \ref{Linalg_Lemma}. Finally, a  simple computation shows
that $A$ is not Codazzi, which is again a contradiction. Thus all
these manifolds cannot be immersed isometrically into the
4-dimensional Euclidean space.\qed
 \subsection*{Acknowledgment} The authors would like to thank Oussama Hijazi for relevant remarks.

\bibliographystyle{amsplain}
\bibliography{immersions}
{\bf Marie-Am\'elie Lawn}\\
Unit\'e de recherche en math\'ematiques\\
Universit\'e de Luxembourg\\
162 A, rue de la fa\"iencerie, L-1511 Luxembourg \\
email: \url{marie-amelie.lawn@uni.lu}\\\\
{\bf Julien Roth}\\
Laboratoire d'Analyse et de Math\'ematiques Appliqu\'ees\\
Universit\'e Paris-Est Marne-la-Vall\'ee\\
Cit\'e Descartes, Champs-sur-Marne, 77454 Marne-la-Vall\'ee Cedex 2, France \\
email: \url{julien.roth@univ-mlv.fr}

\end{document}